\documentclass[runningheads,a4paper]{llncs}

\usepackage{amssymb}
\usepackage{graphicx}
\usepackage{subfigure}
\usepackage{amsmath}

\usepackage{url}
\urldef{\mailsa}\path|{tapanpradhan78,|
\urldef{\mailsb}\path|aurobinda.routray,|
\urldef{\mailsc}\path|bibek.kabi}@gmail.com|
\newcommand{\keywords}[1]{\par\addvspace\baselineskip
\newcommand{\norm}[2]{\| #1 \|}
\noindent\keywordname\enspace\ignorespaces#1}

\begin{document}

\index{Pradhan, Tapan}
\index{Routray, Aurobinda}
\index{Kabi, Bibek}
\mainmatter  

\title{Comparative Evaluation of Symmetric SVD Algorithms for Real-time Face and Eye Tracking}

\titlerunning{Fast SVD}

%
%
\author{Tapan Pradhan%
\and Aurobinda Routray \and\\
Bibek Kabi \\
}
\authorrunning{Tapan Pradhan\and Aurobinda Routray \and\\
Bibek Kabi \\}

\institute{Department of Electrical Engineering, IIT Kharagpur,\\
Kharagpur, India - 721302\\
\mailsa
\mailsb
\mailsc
}
%
\toctitle{Fast SVD}
\tocauthor{Authors' Instructions}
\maketitle
\begin{abstract}
 Computation of singular value decomposition (SVD) has been a topic of concern by many numerical linear algebra researchers. Fast  SVD has been a very effective tool in computer vision in a number of aspects, such as: face recognition, eye tracking etc. At the present state of the art fast and fixed-point power efficient SVD algorithm needs to be developed for real-time embedded computing. The work in this paper is the genesis of an attempt to build an on-board real-time face and eye tracking system for human drivers to detect loss of attention due to drowsiness or fatigue. A major function of this on-board system is quick customization. This is carried out when a new driver comes in. The face and eye images are recorded while instructing the driver for making specific poses. The eigen faces and eigen eyes are generated at several resolution levels and stored in the on-board computer. The discriminating eigen space of face and eyes are determined and stored in the on-board flash memory for detection and tracking of face and eyes and classification of eyes (open or closed) as well. Therefore, fast SVD of image covariance matrix at various levels of resolution needs to be carried out to generate the eigen database. As a preliminary step, we review the existing symmetric SVD algorithms and evaluate their feasibility for such an application. In this article, we compare the performance of (1) Jacobi's, (2) Hestenes', (3) Golub-Kahan, (4) Tridiagonalization and Symmetric QR iteration and (5) Tridiagonalization and Divide and Conquer algorithms. A case study has been demonstrated as an example.
\keywords{Fast SVD, Eigen Space, Face and Eye Detection}
\end{abstract}

\section{Introduction}
SVD is a powerful tool in digital signal and image processing. It states that a matrix can be decomposed as follows,

\begin{equation}
A=U\Sigma V^T
\end{equation}

Where $A_{m\times n}$ is a dense matrix, $U_{m\times m}$ and $V_{n\times n}$ are orthogonal (or unitary) matrices and their columns are called left and right singular vectors respectively. $\Sigma$ is a diagonal matrix and contains all singular values along its diagonal in a non-increasing order. For a symmetric matrix $m=n$ and $U$ and $V$ span the same vector space. Hence computation of either $U$ or $V$ is sufficient. For any dense symmetric matrix $A_{n\times n}$, Eigen Value Decomposition (EVD) is defined as follows

\begin{equation}
A=X \Lambda X^T
\end{equation}

Where $X$ is the matrix of eigenvectors and $\Lambda$ is the diagonal matrix containing eigenvalues along its diagonal.
For symmetric matrices, eigenvalue decompositions and singular value decompositions are closely related as follows \cite{demmel97}:\\
Suppose that $A$ is a symmetric matrix, with eigenvalues $\lambda_i$ and orthonormal eigen vectors $u_i$ so that $A=U \Lambda U^T$ is an eigenvalue decomposition of $A$, with
$\Lambda=diag[\lambda_1 \lambda_2 ...\lambda_n], U=[u_1 u_2 ... u_n]$ and $UU^T=I$. Then an SVD of symmetric matrix $A$ is, $A=U \Sigma V^T$, where diagonal elements of $\Sigma$ i.e. $\sigma_i=abs(\lambda_i)$ and $v_i=sign(\sigma_i).u_i$ where $sign(0)=1$. For symmetric positive definite matrices eigenvalue decomposition (EVD) and SVD leads to the same decomposition. Hence we will use eigenvalues/eigenvectors and singular values/singular vectors interchangeably.\\

However, SVD has been an offline tool for digital signal and image processing applications for decades because of the computation complexity and memory requirement. Due to the increased resources of some of the recently introduced workstations, there are attempts to develop faster versions of SVD algorithm for real-time signal and image processing applications.  The implementation of SVD in embedded platforms like DSPs, ARM and FPGAs is necessary for facilitating efficient real-time image processing. Most of these platforms either have fixed-point processors or CLBs (Configurable Logic Blocks) to make the system cheaper and power efficient. Hence \emph{fast and fixed-point SVD} algorithms are to be developed for such applications. The purpose of this work is to evaluate the existing SVD algorithms for their suitability on embedded platforms.\\

In pattern recognition, eigenspace based method has been proposed for face tracking or face recognition \cite{ahuja2001}-\cite{turk91}. To find the eigenspace, SVD (or eigenvalue decomposition) is used. There are several algorithms for SVD as stated in literature \cite{demmel97}-\cite{datta95}. Jacobi's algorithm is known to be the oldest and slowest algorithm \cite{demmel97}-\cite{parlett98}.For symmetric matrices, though Jacobi's algorithm generates accurate singular values and singular vectors, the time of execution increases with the dimension of the matrices and is only suitable as an offline tool.Two-sided and one-sided variants for Jacobi's algorithm are stated in literature. Hestenes' algorithm is a variant of one-sided Jacobi's algorithm and is discussed in \cite{hestenes58}-\cite{svensson89}. Being a one-sided version the time of computation is lesser than two-sided Jacobi's algorithm. However, as the iteration is applied on the whole process this algorithm is also not suitable for online applications. Golub and Kahan proposed a two-step algorithm \cite{golubloan97},\cite{golubkahan65}-\cite{demmelkahan90} for computation of SVD. In the first step, a dense symmetric matrix is converted to a bidiagonal matrix, which is eventually converted to a diagonal matrix using implicit QR iteration in the second step. Because the second phase of Golub-kahan algorithm is iterative in nature, it is much faster than Jacobi's or Hestenes' algorithm. A similar two step algorithm has been proposed for SVD, where a dense symmetric matrix is reduced to tridiagonal matrix and then an implicit symmetric QR iteration is applied to reduce the symmetric tridiagonal matrix into a diagonal matrix. This algorithm is found to be faster and competitive with Golub-Kahan algorithm when a combination of QR and QL algorithm is used \cite{demmel97}-\cite{datta95}. Still these algorithms could not satisfy real-time constraints as required by the signal and image processing platforms. Hence a  Divide and Conquer algorithm was proposed by J.J.M. Cuppen based on a rank-one modification by Bunch, Nielsen and Sorensen. This is the fastest algorithm till date when a complete eigensytem of a symmetric tridiagonal matrix is required \cite{demmel97}. A variant of the said algorithm by Gu and Eisenstat has been implemented  in LAPACK routine for matrices with dimension larger than 25 \cite{bns78}-\cite{cuppen81}.\\

 Faster performance is achieved when floating-point SVD is converted to fixed-point format and implemented in fixed-point platform \cite{zoran2007}. Fast and Fixed-point SVD algorithm is also useful for reducing silicon area and power consumption in embedded platforms \cite{kim98}-\cite{keding98}. For Digital Signal Processing applications, attempts have been made to implement SVD algorithm using multiprocessor arrays \cite{brent85} and CORDIC (COordinate Rotation DIgital Computer) based reconfigurable systems \cite{wang2008}-\cite{szecowka10}.\\

\section{Existing Algorithms for SVD of Symmetric Matrices }
In this section, we analyse different SVD algorithms for symmetric matrices along with their complexity, advantages and disadvantages.
\begin{enumerate}
\item Jacobi's Algorithm \cite{parlett98}
\item Hestenes' Algorithm \cite{svensson89}
\item Golub-Kahan Algorithm \cite{golubkahan65}-\cite{demmelkahan90}
\item Tridiagonalization and Symmetric QR iteration \cite{demmel97}-\cite{datta95}
\item Tridiagonalization and Divide and Conquer Algorithm \cite{demmel97}
\item Bisection and Inverse Iteration \cite{demmel97}
\end{enumerate}

\begin{table}
\begin{center}
\caption{Computational complexities of SVD algorithms}
\begin{tabular}{|l|r|}
\hline
Algorithms & Time Complexity\\
\hline
Jacobi's Algorithm & $O(n^3)$\\
Hestenes Method & $O(n^3)$\\
\hline
Golub-Kahan Algorithm  & ($\Sigma[(\frac{8}{3}n^3 + O(n^2))]$ and \\
(Bidiagonalization + Implicit QR Iteration) & $U,V[4n^3+O(n^3)]$)\\
Tridiagonalization + Symmetric QR Iteration & $8\frac{2}{3}n^3 + O(n^2)$\\
Tridiagonalization + Divide and Conquer Method & $\frac{8}{3}n^3 + O(n^2)$\\
\hline
\end{tabular}
\end{center}
\label{tab:complexity}
\end{table}
\subsection{Jacobi's Algorithm }
Jacobi's algorithm is the oldest and slowest available method for computing SVD by implicitly applying iteration on a dense symmetric matrix $A$.
The method is more or less similar to the method for eigenvalue decomposition of symmetric matrices. Jacobi's method computes SVD of a symmetric matrix A with higher relative accuracy when A can be written in the form A=DX, where D is diagonal and X is well conditioned \cite{demmel97}. Thus
\begin{equation}
J^TAJ=\Sigma
\end{equation}
In each step we compute a Jacobi rotation with $J$ and update $A$ to $J^TAJ$, where $J$ is chosen in such a way that two off-diagonal entries of a $2\times2$ matrix of $A$ are set to zero in $J^TAJ$ . This is called two-sided or classical Jacobi method. Again, this can be achieved by forming $G=A^TA$ and performing iteration over $G$ instead of $A$. The eigenvalues of $G$ are the square of the singular values of A. Since $J^TGJ=J^TA^TAJ=(AJ)^T(AJ)$, we can obtain $\Sigma$ by merely computing $AJ$ . This is termed as one-sided Jacobi rotation. Though Jacobi's algorithm for calculating SVD is the slowest among all presently available algorithms, it produces relatively accurate singular values and singular vectors. Jacobi's algorithm can be implemented in parallel as the individual steps are not interdependent. Parallel Jacobi's algorithm is discussed in \cite{zhoubrent1}.

\begin{table}
\begin{center}
\caption{Two-sided Jacobi's algorithm for SVD}
\begin{tabular}{|l|l|}
\hline
$Q=I$ ($I =$ identity matrix) & \dots\\
 repeat & s=c.t\\
for $j=1:(n-1)$ & $J(j,j)=c, J(k,k)=c,$ $J(j,k)=s, J(k,j)=-s$\\
for $k=(j+1):n$ & $A=J^TAJ$\\
if $|A(j,k)|$ is not too small & $Q=Q.J$\\
$\xi=\frac{(A(k,k)-A(j,j))}{2A(j,k)}$ & endif\\
if $\xi=0$ & endfor\\
$t=\frac{1}{(|\xi|+\sqrt{1+\xi^2}}$ & endfor\\
else & $\sigma_k=|A(k,k)|$\\
$t=\frac{sign(\xi)}{(|\xi|+\sqrt{1+\xi^2}}$ & $U=Q$\\
endif & $v_k=sign(A(k,k)).u_k$\\
 $c=\frac{1}{\sqrt{1+t^2}}$ & $V=[v_1 v_2 \dots v_n]$\\
contd \dots & \\
\hline
\end{tabular}
\end{center}
\label{tab:jacobi2s}
\end{table}
\subsection{Hestenes' Algorithm }
Hestenes' algorithm is based on one-sided Jacobi's algorithm. A symmetric matrix of size $n\times n$ is used to generate an orthogonal rotation matrix $V$ so that the transformed matrix $A'=AV=W$ has orthogonal columns. Now if we normalize each non-null column of matrix $W$ to unity, we get the relation $W=AV=U\Sigma=U\cdot diag(\sigma_1,\sigma_2,...,\sigma_n)$,
$A=U\Sigma V^T, V=\prod_{i=1}^n J_i$ and singular values of $A$ are $\sigma_i= {\parallel {a'_i}\parallel}_2$, where $A'=[a'_1a'_2\dots a'_n]$.\\
\begin{table}
\begin{center}
\caption{Hestenes' algorithm for SVD}
\begin{tabular}{|l|l|}
\hline
$V=I$ & \dots\\
repeat until convergence & $c=\frac{1}{\sqrt{1+t^2}}, s=t.c$\\
for $p=1:(n-1)$ & end\\
for $q=(p+1):n$ & J=I\\ 
$\alpha={A(:,p)^T}A(:,p)$, $\beta={A(:,q)^T}A(:,q)$ & $J(p,p)=c, J(q,q)=c, J(p,q)=s, J(q,p)=-s$\\
$\gamma={A(:,p)^T}A(:,q)$ & $A=A.J$, $V=V.J$\\
if $\gamma=0$ & endfor\\
$c=1; s=0$ & endfor\\
else & for $k=1:n$\\
$\zeta=\frac{(\beta-\alpha)}{2\gamma}$ & $\sigma_i=||A(:,k)||$\\
 $t=\frac{sign(\zeta)}{(|\zeta|+\sqrt{1+\zeta^2})}$ & end\\
contd \dots & $U=AV\Sigma^{-1}$\\
\hline
\end{tabular}
\end{center}
\label{tab:hestenes}
\end{table}

\subsection{Golub-Kahan Algorithm }
This algorithm can be segmented in two phases,\\
i)	In the first phase a dense symmetric matrix is converted to a bi-diagonal matrix by orthogonal similarity transformation.\\
ii)	This bi-diagonal matrix   is then converted to diagonal matrix using Implicit QR iteration.\\

The standard algorithm may compute singular values with poor relative accuracy. However, with modified algorithm by Demmel-Kahan smaller singular values may be computed with high relative accuracy \cite{demmelkahan90}.
\begin{table}
\begin{center}
\caption{Golub-Kahan algorithm for SVD}
\begin{tabular}{|l|l|}
\hline
function $[u,\sigma]$=\textbf{houszero}(x) & function $[c,s,r]$=\textbf{rot}(f,g)\\
$m=max(|x_i|), i=1,2,\dots,n$ & if $f=0$ then $c=0$; $s=1$; $r=g$\\
$u_i=\frac{x_i}{m}, i=1,2,\dots,n$ & elseif ($|f|>|g|$) then\\
$sign(0)=1$ & $t=\frac{g}{f}; tt=\sqrt{1+t^2}; c=\frac{1}{tt}; s=t.c; r=tt.f$\\
$\sigma=sign(u_1)\sqrt{{u_1}^2+{u_2}^2+\dots+{u_n}^2}$ & else\\
$u_1=u_1+\sigma$ & $t=\frac{f}{g}; tt=\sqrt{1+t^2}; s=\frac{1}{tt}; c=t.s; r=tt.g$\\
$\sigma=-m.\sigma$ & endif\\
end \textbf{houszero} & end \textbf{rot}\\
\hline
\textbf{Golub-Kahan algorithm for SVD} & \dots\\
$U_1=I, V_1=I$ &  repeat \\
for $k=1:(n-1)$ &  for $i=1:(n-1)$\\
$[u,\sigma]=\textbf{houszero}(A(k:m,k)$ & $U_2=I$\\
$H1=I-2\frac{uu^T}{u^Tu}$ &  $[c,s,r]=$\textbf{rot}($A(i,i),A(i,i+1)$)\\
$P_1=I$ &  $Q=I$\\
$P_1(k:m,k:n)=H1$ & $Q(i:i+1,i:i+1)=
\begin{pmatrix}
c & s\\
-s & c\\
\end{pmatrix}$\\
$A(k:m,k:n)=H1.A(k:m,k:n)$ & $A=AQ^T$ \\
$U_1=U_1 . P_1$ & $V_2=V_2.Q^T$\\
if $k \leq (n-2)$ & $[c,s,r]=$\textbf{rot}($A(i,i),A(i+1,i)$)\\
$[v,\sigma]=$\textbf{houszero}($A(k,k+1:n)^T$) & $Q=I$\\
$H2=I-2\frac{vv^T}{v^Tv}$ & $Q(i:i+1,i:i+1)=
\begin{pmatrix}
c & s\\
-s & c\\
\end{pmatrix}$\\
$P_2=I$ & $A=QA$\\
$P_2(k:m-1,k:n-1)=H2$ & $U_2=U_2.Q$\\
$A(k:m,k+1:n)=A(k:m,k+1:n).H2$ & endfor\\
$V_1=V_1.P_2$ & $\Sigma=abs(A)$\\
endif & $U^T=U_1.U_2$\\
endfor & $V=V_1.V_2$\\
contd \dots & \\
\hline
\end{tabular}
\end{center}
\label{tab:golubkahan}
\end{table}
This algorithm has two steps as shown in equation (\ref{eq:gkalgo}).
\begin{equation}\label{eq:gkalgo}
\begin{pmatrix}
* & * & * & * \\
* & * & * & * \\
* & * & * & * \\
* & * & * & * \\
\end{pmatrix}
\xrightarrow[\text{STEP I}]{U_1A{V_1}^T=B}
\begin{pmatrix}
* & * & 0 & 0 \\
0 & * & * & 0 \\
0 & 0 & * & * \\
0 & 0 & 0 & * \\
\end{pmatrix}
\xrightarrow[\text{STEP II}]{U_2B{V_2}^T=\Sigma}
\begin{pmatrix}
* & 0 & 0 & 0 \\
0 & * & 0 & 0 \\
0 & 0 & * & 0 \\
0 & 0 & 0 & * \\
\end{pmatrix}
\end{equation}
\begin{center}
\end{center}

Hence, ${U_2}^T B V_2={U_2}^TU{_1}^TA V_1 V_2={(U_1U_2)}^T A{V_1 V_2}=U^T A V=\Sigma$ and $\Sigma$ is the diagonal matrix containing singular values of the symmetric matrix $A$. In bi-diagonalisation process we have used householder method to make the elements of a column or row zero. Another optimized variant of this algorithm is also available and is known as Golub-Kahan-Chan algorithm.
\subsection{Tridiagonalization and Symmetric QR iteration }
In this method, the dense symmetric matrix is first converted to symmetric tridiagonal matrix in finite number of steps and the tridiagonal matrix is eventually converted to diagonal form by symmetric QR iteration. When a combination of QL and QR algorithm instead of QR is used, a stable variant is obtained. The steps of this algorithm is stated in equation (\ref{eq:symqr}).\\

\textbf{Method of tridiagonalization}\\
Householder reflection can be used for tridiagonalization. In this method, one row or column is picked up, and householder matrix is reconstructed to make all the elements of the row or column zero except the first one \cite{golubloan97}.\\

Hessenberg reduction [\cite{datta95}] is a process by which a dense matrix is converted to an upper or lower Hessenberg matrix by orthogonal similarity transformation and thus does not change the eigenvalues or singular values. For a dense symmetric matrix, this Hessenberg reduction produces an upper as well as a lower Hessenberg matrix or a symmetric tridiagonal matrix which is used in symmetric QR iteration process to produce a diagonal matrix.In this process, the eigenvalues or singular values remain same as that of the original dense symmetric matrix. Householder or Givens method may be applied for Hessenberg reduction. We have used Givens method to produce symmetric tridiagonal matrix.\\

Lanczos method is useful to transform a symmetric dense matrix to a symmetric tridiagonal matrix. This method suffers from loss of orthogonality among Lanczos vectors with increased number of iterations. For this reason reorthogonalization is required to obtain proper orthogonal vectors.\\

\emph{Symmetric QL and QR iteration} :\\
This is an iterative process and has a complexity of $O(n^2)$. After reducing the dense symmetric matrix into its tridiagonal form by orthogonal similarity transformation, symmetric QL or QR iteration is performed depending on which of the first and last diagonal elements of the symmetric tridiagonal matrix is larger. If the first diagonal entry is larger than the last one, QR is called upon to perform, else QL is called.
\begin{equation}\label{eq:symqr}
\begin{pmatrix}
* & * & * & * \\
* & * & * & * \\
* & * & * & * \\
* & * & * & * \\
\end{pmatrix}
\xrightarrow[\text{STEP I}]{U_1A{V_1}^T=T}
\begin{pmatrix}
* & * & 0 & 0 \\
* & * & * & 0 \\
0 & * & * & * \\
0 & 0 & * & * \\
\end{pmatrix}
\xrightarrow[\text{STEP II}]{U_2T{V_2}^T=\Sigma}
\begin{pmatrix}
* & 0 & 0 & 0 \\
0 & * & 0 & 0 \\
0 & 0 & * & 0 \\
0 & 0 & 0 & * \\
\end{pmatrix}
\end{equation}
\begin{center}
\end{center}
\begin{table}
\begin{center}
\caption{Tridiagonalization and Symmetric QR Iteration}
\begin{tabular}{|l|l|}
\hline
function $[c,s]=$\textbf{givens}($a,b$) & $P_1(k+1:n,k+1:n)=H$\\
if $b=0$ & $A=P_1AP_1$\\
$c=1; s=0$ & $U_1=U_1.P_1$\\
else & end\\
if $|b|>|a|$ & repeat\\
$\tau=-\frac{a}{b}$ & $d=\frac{(t_{n-1,n-1}-t_{n,n})}{2}$\\
$s=\frac{1}{\sqrt{1+\tau^2}};c=s.\tau$ & $\mu=t_{n,n}-\frac{t_{n,n-1}^2}{(d + sign(d)\sqrt{d^2+t_{n,n-1}^2})}$\\
else & $x=t_{11}-\mu; z=t_{21}$\\
$\tau=-\frac{a}{b}$ & for $k=1:(n-1)$\\
$c=\frac{1}{\sqrt{1+\tau^2}}; s=\frac{c}{\tau}$ & $[c,s]$=\textbf{givens}$(x,z)$\\
endif & $T=(G_k^T)TG_k$, where $G_k = G(k,k+1,\theta)$\\
.............................................. & $U_1 = U_1G_k$\\
$U_1=I$ & if $k < (n-1)$\\
for $k=1:(n-2)$ & $x=t_{k+1,k}$\\
$[u,\sigma]$ =\textbf{houszero}($x$) & $z = t_{k+2,k}$\\
$H=I_1-2\frac{(uu^T)}{(u^Tu)}$ & endif\\
$P_1=I$ & endfor\\
$contd..$ & $\Sigma=abs(T), U=U_1$\\
\hline
\end{tabular}
\end{center}
\label{tab:tridiagqr}
\end{table}
\subsection{Tridiagonalization and Divide and Conquer Algorithm }
The method was first proposed by J.J.M. Cuppen \cite{cuppen81} based on a rank-one modification by Bunch, Nielsen and Sorensen \cite{bns78}. However, the algorithm became popular after a stable variant for finding singular vectors or eigenvectors was found out in 1990 by Gu and Eisenstat \cite{gueisenstat92}.
This algorithm is the fastest SVD algorithm available till date. A dense symmetric matrix is first converted to symmetric tridiagonal matrix \cite{demmel97}, \cite{bns78}-\cite{cuppen81}.
Then the symmetric tridiagonal matrix is divided into two parts by rank one update and again each of the smaller matrices is divided till sufficiently smaller (matrix dimension = 25) matrices are formed. Then QR and QL iteration may be applied to find the SVD of the smaller matrices and using rank one update smaller solutions are combined together to form the complete SVD of the symmetric tridiagonal matrix. With a combination of previously stated steps, a complete eigensystem of a dense symmetric matrix can be found out.\\

Two major parts for finding the eigensystem of a symmetric tridiagonal matrix are \emph{divide} and \emph{conquer}. It works by breaking down a problem into two or more subproblems of the same type until the subproblem becomes simple enough to be solved directly. The solutions to the subproblems are then combined to generate a solution to the original problem. The most significant part is the solution of secular equation. The solution of the secular equation involves function approximation for finding the desired roots. For matrices with dimension greater than 25 this is the fastest method for finding the complete eigensystem of a symmetric tridiagonal matrix till date \cite{demmel97}.
\subsubsection{Structure of the Algorithm \cite{demmel97}}
This includes dividing the symmetric tridiagonal matrix into two parts by removing the subdiagonal entry by rank one modification. Now with the known eigensystem of the two new symmetric tridiagonal matrices the secular equation is constructed. Solution of secular equation gives the eigenvalues of the original matrix from which the eigen vectors can also be computed. A dense symmetric matrix is converted to a symmetric tridiagonal matrix and then the divide and conquer algorithm is applied  as \cite{gueisenstat92}.
\begin{enumerate}
\item Step 1 - Dividing the symmetric tridiagonal matrix into smaller matrices by rank one modification.
\item Step 2 - Sorting eigenvalues and eigenvectors obtained  from smaller matrices in an increasing order using permutation matrix.
\item Step 3 - Deflation (when updation of eigenvalues and eigen vectors is not required) due to smaller coefficient or smaller components in the vector used for rank one modification and also for repeated eigenvalues are considered.
\item Step 4 - Formation and solution of secular equation using non-deflated eigenvalues.
\item Step 5 - Combining the solution obtained from secular equation solver and deflated eigenvalues to obtain the complete eigensystem of the symmetric tridiagonal matrix.
\end{enumerate}
\emph{Forming the symmetric tridiagonal matrix from a dense symmetric matrix}\\
  \begin{equation}
U^TAU=T
\end{equation}
$A$ - Dense symmetric matrix\\
$T$ - Symmetric tridiagonal matrix\\
$U$ - Eigenvector matrix obtained from tridiagonalization process\\

\textbf{Divide and Conquer Algorithm}\\
\emph{Step 1 Rank one modification of the symmetric tridiagonal matrix}\\
\begin{equation}
T =
\begin{pmatrix}
T_1 & 0\\
0 & T_2
\end{pmatrix}
+ b_m vv^T = \\
\begin{pmatrix}
Q_1D_1{Q_1}^T & 0\\
0 & Q_2D_2{Q_2}^T
\end{pmatrix}
+ b_m vv^T
\end{equation}
\begin{equation}
T =
\begin{pmatrix}
Q_1 & 0\\
0 & Q_2
\end{pmatrix}
\begin{pmatrix}
\begin{pmatrix}
D_1 & 0\\
0 & D_2
\end{pmatrix}
+ b_m uu^T
\end{pmatrix}
\begin{pmatrix}
{Q_1}^T & 0\\
0 & {Q_2}^T
\end{pmatrix}
\end{equation}
Let, $D=
\begin{pmatrix}
D_1 & 0\\
0 & D_2
\end{pmatrix}
, b_m=\rho
$\\
Where, $T_1, T_2$ are smaller tridiagonal matrices.\\
$D$ is a diagonal matrix containing the eigenvalues of the smaller tridiagonal matrices.\\
$Q_1, Q_2$ are eigenvector matrices after eigendecomposition of smaller tridiagonal matrices.\\
$u=
\begin{pmatrix}
{Q_1}^T & 0\\
0 & {Q_2}^T
\end{pmatrix}
v=
\begin{pmatrix}
last & column & of & {Q_1}^T\\
 first & column & of & {Q_2}^T
\end{pmatrix}$\\
\emph{Step 2 Sorting using permutation matrix}\\
Following equation is used for sorting process to arrange the eigenvalues and eigenvectors in an increasing order.
\begin{equation}\label{eq:permut}
  D + \rho u u^T=P^T P(D + \rho u u^T)P^T P
  \end{equation}
Where P is the permutation matrix.\\
\emph{Step 3 Reducing the computation using deflation}\\
Deflation occurs when \\
1. $\rho$ is very small, 2. smaller weights ($u_i$) due to smaller components in the vector used for rank one modification and 3. multiple eigenvalues \cite{bns78}. When deflation occurs eigenvalues and eigenvectors do not require updation. Hence a saving in computation is achieved.\\
\emph{Step 4 Formation and solution of secular equation}\\
 $|((D + \rho u u^T) - \lambda I)|=0$\\
 $|((D-\lambda I)(I + \rho (D-\lambda I)^{-1} u u^T)|=0$\\
Since, $|(D- \lambda I)|\not = 0$, $|(I + \rho (D-\lambda I)^{-1} u u^T)|=0$\\
Now,
\begin{equation} \label{eq:seceq}
|((I + \rho (D-\lambda I)^{-1} u u^T)|= 1 + \rho u^T (D-\lambda I)^{-1} u= 1+ \rho{\sum_{i=1}^n \frac{{u_i}^2}{d_i-\lambda}}
\end{equation}
\begin{equation}
\begin{pmatrix}
D_1 & 0\\
0 & D_2
\end{pmatrix}
+ b_m uu^T=[R][\Lambda][R^T]
\begin{pmatrix}
obtained by & solving & secular & equation
\end{pmatrix}
\end{equation}
Eigenvalues and eigenvectors are obtained by solving a secular equation like equation (\ref{eq:seceq}) \cite{bns78}-\cite{cuppen81}.\\
\emph{Step 5 Combining the solutions from Step IV and V the complete eigensystem of the symmetric tridiagonal matrix is obtained}\\
\begin{equation}
T=[Q]
\begin{pmatrix}
\begin{pmatrix}
D_1 & 0\\
0 & D_2
\end{pmatrix}
+ b_m uu^T
\end{pmatrix}
{[Q]}^T=
[Q][R][\Lambda]{[R]}^T{[Q]}^T=[X][\Lambda]{[X]}^T
\end{equation}
Where, $[X]=[Q][R]$.\\
Hence,
\begin{equation}
A=[U][T]{[U]}^T=[U][X][\Lambda]{[X]}^T{[U]}^T=[V][\Lambda]{[V]}^T
\end{equation}
There are some issues related to divide and conquer algorithm. They are discussed below.
\subsubsection{Sorting with Permutation Matrices \cite{jrutterlawn69}}
 If $d_1<d_2\dots<d_n$ then the sequence of eigenvalues obtained will be $\lambda_1<\lambda_2<\dots <\lambda_n$. However, we may not come across a diagonal matrix with sorted diagonal elements after rank one modification and QR or QL iteration. Hence, we need to apply permutation to sort them in ascending order using permutation matrix.\\
 Let, $D=
 \begin{pmatrix}
 13.1247 & 0 & 0 & 0\\
 0 & 201.9311 & 0 & 0\\
 0 & 0 & 0.0693 & 0\\
 0 & 0 & 0 & 26.7189
 \end{pmatrix}
 $\\
 Let, $ u=
 \begin{pmatrix}
 -0.5421\\
 -0.4540\\
 0.2128\\
 -0.6743
 \end{pmatrix}
 $\\
 Now applying permutation using permutation matrix,
 $P=
 \begin{pmatrix}
  0 & 0 & 1 & 0\\
  1 & 0 & 0 & 0\\
  0 & 0 & 0 & 1\\
  0 & 1 & 0 & 0
  \end{pmatrix}
  $
  and using equation (\ref{eq:permut}) , we obtain
  $D(sorted)=
  \begin{pmatrix}
  0.0693 & 0 & 0 & 0\\
  0 & 13.1247 & 0 & 0\\
  0 & 0 & 16.7189 & 0\\
  0 & 0 & 0 & 201.9311
  \end{pmatrix}
  $
  and
  u(modified)=
  $
  \begin{pmatrix}
  0.2128\\
  -0.5421\\
  -0.6743\\
  -0.4540
  \end{pmatrix}
  .$
 \subsubsection{Solving Secular Equation}
 The roots of the secular equation $1+ \rho{\sum_{i=1}^n \frac{{u_i}^2}{d_i-\lambda}}$ are the eigenvalues of the original matrix. $\rho$ (rho)  is the sub diagonal entry which creates the rank one modification. ${u_i}^2$ is the weight over the pole. The secular equation has poles at the eigenvalues of D and zeros at the eigenvalues of $D + \rho u u^T$.\\
 According to interlacing property :
\begin{enumerate}
\item If rho is greater than zero, the roots lie in such a manner that : $d_1<\lambda_1<d_2<\lambda_2<\dots<d_n<\lambda_n$
\item If rho is less than zero then :  $\lambda_1<d_1<\lambda_2<d_2<\dots<\lambda_n<d_n$
\end{enumerate}

Assuming rho is greater than zero for $i<n$, the roots lie in between $d_i$ and $d_{i+1}$, but for $i=n$, the root lies in a manner that $d_n<\lambda_n<d_n + \rho uu^T$ \cite{bns78}.
For the given matrix
$
\begin{pmatrix}
16.7118 & 10.7270\\
10.7270 & 34.2341
\end{pmatrix}$
, where rho (10.7270) is greater than zero, the nature of the roots are examined below

\begin{figure}
\begin{center}
\fbox{
\includegraphics[width=0.8\textwidth]{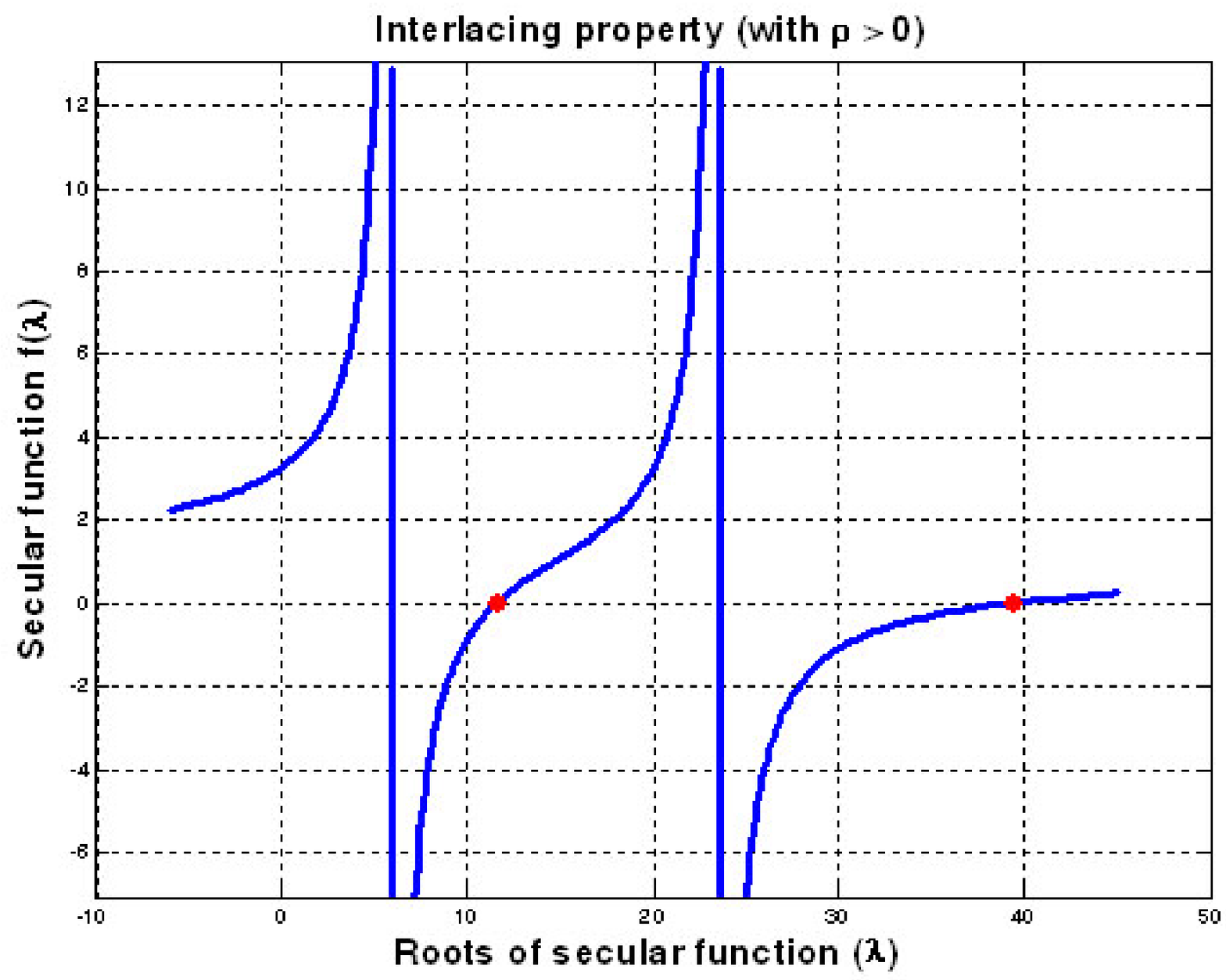}
}
\caption{Roots showing interlacing property with $rho>0$} \label{fig:cfig1}
\end{center}
\end{figure}
\begin{figure}
\begin{center}
\fbox{
\includegraphics[width=0.8\textwidth]{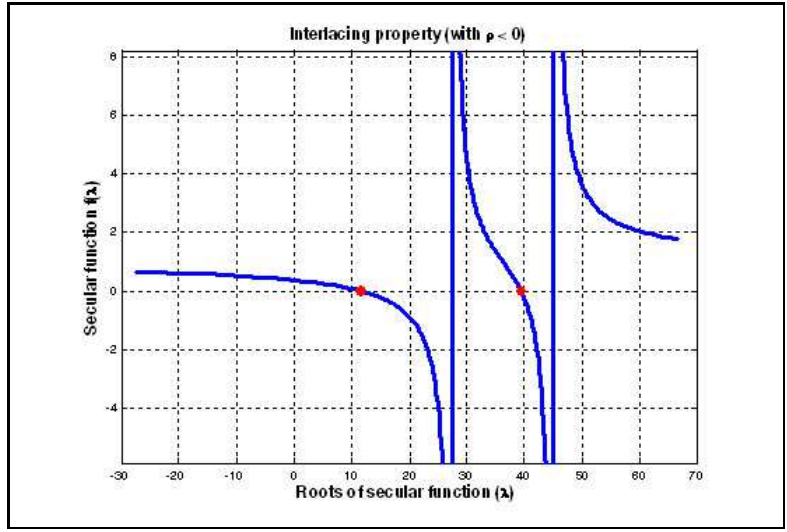}
}
\caption{Roots showing interlacing property with $rho<0$} \label{fig:cfig2}
\end{center}
\end{figure}

The eigenvalues of the previously stated matrix are
$
\begin{pmatrix}
11.6228 & 0\\
0 & 39.3231
\end{pmatrix}$
, where the diagonal matrix $D$ is
$
\begin{pmatrix}
5.9848 & 0\\
0 & 23.5071
\end{pmatrix}.$
In Figure ~\ref{fig:cfig1} we can see bold  lines and  vertical dotted lines. The points, where the bold lines cross the real axis at zero, are the roots for the corresponding matrix. The vertical dotted lines represent the diagonal elements after the rank one modification and eigendecomposition of smaller matrices.\\
For a matrix
$
\begin{pmatrix}
16.7118 & -10.7270\\
-10.7270 & 34.2341
\end{pmatrix}$
with rho (-10.7270) less than zero, the roots appear as shown in Figure ~\ref{fig:cfig2}.

The eigenvalues of the above said matrix with $rho<0$ are
$
\begin{pmatrix}
11.6228 & 0\\
0 & 39.3231
\end{pmatrix}$
, where the diagonal matrix $D$ is
 $
\begin{pmatrix}
27.4388 & 0\\
0 & 44.9611
\end{pmatrix}.$
The interlacing property is clearly found from Figures ~\ref{fig:cfig1} and ~\ref{fig:cfig2}.
\subsubsection{Methods of Solving Secular Equation \cite{li94}}
We assume that $ \rho $ (rho) is greater than zero. The solution of secular equation is based on function approximation. So we may think of using Newton's iterative procedure to solve the equation, but Newton's method is based on local linear interpolation. At some points of initial guess $(\lambda_0)$ for the desired root, the linear approximation would be horizontal and the next approximation would be a large negative number which is not useful [5]. This happens particularly when the weights are small. Therefore, we go for rational osculatory interpolation. Rational osculatory interpolation of secular function is the combination of two kinds of rational functions. If the secular function is $f(\lambda)$, the two rational functions are respectively $\Phi(\lambda)$ and $\Psi(\lambda)$. $\Phi(\lambda)$ is the sum of positive terms and $\Psi(\lambda)$ is the sum of negative terms. The inequality among them is $-\infty<\Psi(\lambda)<0<\Phi(\lambda)<+\infty$. In the following section, we will explore a set of schemes for solving secular equation. The first scheme is named as \emph{approaching from left} because the algorithm will produce a sequence of monotonically increasing approximations to desired root provided the initial guess is less than the desired root. Starting with the initial guess which is less than the desired root, the scheme will yield a sequence of approximations converging monotonically upwards to the desired root. On the other hand, if the initial guess exceeds the required root so much that the next approximation does not converge, the scheme fails to get the desired root. \emph{Approaching from left} scheme proves better for $i<n$, but for $i=n$ according to the interlacing property the  root lies in the range $d_n< \lambda_n< d_n + \rho u u^T$. Here the initial guess exceeds the desired root and the second scheme i.e.\emph{approaching from right} will yield a sequence of approximations converging monotonically downwards to desired root. A \emph{middle way method} has been adopted by combining both the above schemes. This scheme takes both nearby poles into consideration, while the other two schemes take each of the nearby poles into consideration. Here the scheme handles two parts. The first part is, if the combination of two rational functions is greater than zero, the root is closer to $d_i$ and if the function is less than zero, the root is closer to $d_{i+1}$.
Interpolation of both $\Psi(\lambda)$ and $\Phi(\lambda)$ is carried out by taking both nearby poles into consideration.
Following \cite{demmel97},
 \begin{equation}
h(\lambda)=\frac{c_1}{d_i-\lambda} + \frac{c_2}{d_{i+1}-\lambda} + c_3
\end{equation}
\begin{center}
$
f(\lambda)=1 + \sum_{k=1}^i\frac{{u_k}^2}{d_k-\lambda} + \sum_{k=i+1}^n\frac{{u_k}^2}{d_{k}-\lambda}
$
\end{center}
\begin{equation}
= 1 + \Phi(\lambda) + \Psi(\lambda)
\end{equation}
Where coefficients $c_1,c_2$ and $c_3$ are computed by interpolating $f(\lambda)$.
In \cite{li94} the above mentioned secular function has been represented similarly using the system's matrix as:
\begin{equation}
D + \frac{1}{\rho} u u^T
\end{equation}
\begin{equation}
h(\lambda)=\rho + r + R + \frac{s}{d_k-\lambda} + \frac{S}{d_{k+1}-\lambda}
\end{equation}
Where $r+ \frac{s}{d_k-\lambda}$ is approximated to $\Psi(\lambda)$ and $R + \frac{S}{d_{k+1}-\lambda}$ is approximated to $\Phi(\lambda)$.
The initial guess $(\lambda_0)$ plays a very significant role in finding solution for secular equation. For $k<n$ the initial guess is calculated as $\frac{d_k + d_{k+1}}{2}$ and for $k = n$ the initial guess is $d_n + \rho u u^T$. So with $\lambda_0$, we calculate the value of $h(\lambda)$,  and if the rational function is greater than zero, the root is closer to $d_i$ and if it is less than zero, the root is closer to $d_{i+1}$. So according to the placement of the root, if it is close to $d_i$ then we shift each of $\lambda_0$, $d_i$ and $d_{i+1}$ by subtracting $d_i$ from them. With the new value of the initial guess $(\lambda_1)$, we solve the quadratic equation formed by the two rational functions, and $\eta$ (iterative correction ) is found. The desired root is found by computing $(\lambda_1 + \eta)$. We know that the weights over the poles play a major role in the solution of secular equation. One of the circumstances is associated with weights. In order to overcome this, fixed weight method was implemented with the structure same as middle way method. By combining \emph{middle way method} and \emph{fixed weight method}, a hybrid scheme was designed to make iteration faster. In this scheme, the function interpolation was carried out with $d_k$, $d_{k+1}$  and $d_{k-1}$. The rational function for hybrid scheme is given below.
\begin{equation}
h(\lambda)=c + \frac{s}{d_{k-1}-\lambda} + \frac{{u_k}^2}{d_k - \lambda}+ \frac{S}{d_{k+1}-\lambda}
\end{equation}
Where $c=\rho + r + R$.
Once the eigenvalues are obtained, eigenvectors can be computed using the following equation
\begin{equation} \label{eq:xvec}
x=\frac{(\lambda I -D)^{-1} u}{\parallel{(\lambda I -D)}^{-1} u \parallel}
\end{equation}
Unfortunately, the above equation does not produce accurate eigenvectors [17,20] and the following equation is used for modifying $u$ vectors where $u$ vectors are updated and  used in equation (\ref{eq:xvec})
 \begin{equation}
 u_k=\sqrt\frac{\prod_{j=1}^{k-1}(d_k - \lambda_j) \prod_{j=k}^n (\lambda_j - d_k)}{\rho \prod_{j=1}^{k-1}(d_k - d_j) \prod_{j=k+1}^n (d_j - d_k)}
 \end{equation}
In this way, the eigenvectors computed (using equation \ref{eq:xvec}) are relatively accurate \cite{gueisenstat92}.\\

Apart from the above SVD algorithms, \emph{Bisection with Inverse Iteration} can also be used. This is specially used when singular values or eigenvalues are required within a specified range. This algorithm suffers from loss of orthogonality among the computed singular vectors or eigenvectors. Hence forced re-orthogonalization is necessary.
\section{Case Study : Real-time Face and Eye Tracking}
\begin{figure}
\centering
\fbox{
\includegraphics[width=5in]{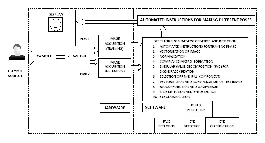}
}
\caption{Schematic diagram of the proposed online training based face and eye tracking system.} \label{fig:dblock}
\end{figure}

 A fast SVD based face tracking scheme based on eigen space is shown in Figure ~\ref{fig:dblock}. The scheme is composed of two sections - hardware and software. Fast and fixed-point SVD is required in this scheme for rapid customization of the system for a particular human subject. There are two paths: path 1 is for training and path 2 is for tracking purposes. The subject is asked to put its face within a box area (frame) in the image. The subject is then required to follow some instructions which direct him to make different poses with different face orientations. The training set is prepared online from extracted face images with different face orientations (front, front up, front down, looking left, looking right etc.) and it is customized for the person sitting in front of the camera. The covariance matrix is constructed from extracted face images with different poses at a frame rate of 30 fps. Fast SVD is performed on the covariance matrix to prepare the eigen space or feature space of face or eye in real-time. Now face and eye detection and tracking is carried out by projecting a block image from the incoming frame (from web camera) onto the feature vector space and comparing the reconstructed image with a standard face or eye image. Face and eyes are detected at approximately 10 fps in the laboratory environment.

\begin{figure}
\begin{center}
\fbox{
\includegraphics[width=0.8\textwidth]{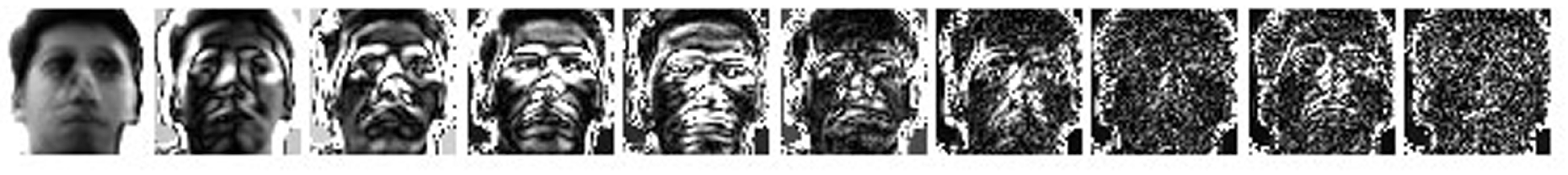}
}
\fbox{
\includegraphics[width=0.8\textwidth]{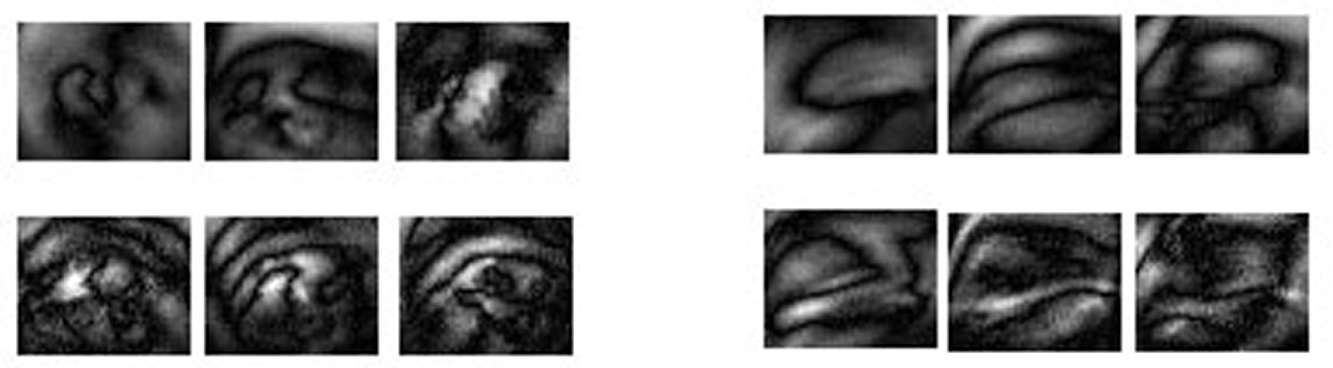}
}

\fbox{
\includegraphics[width=0.8\textwidth]{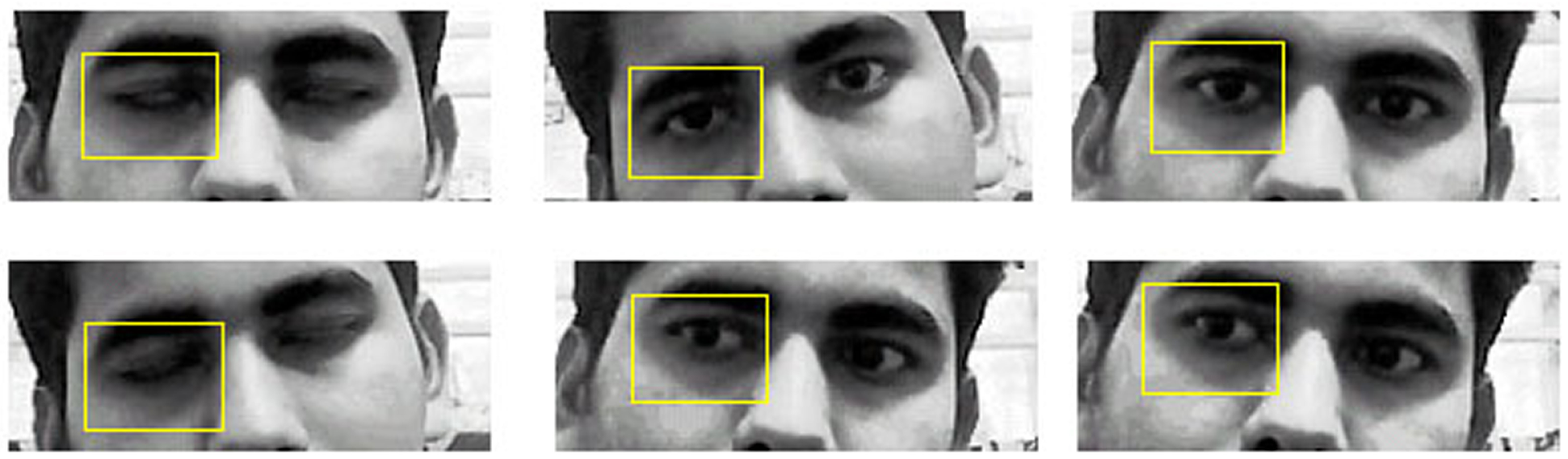}
}
\caption{EigenFaces (top), Eigen Eyes (for open and closed eyes) and Eye Detection using Fast SVD based algorithm (bottom).} \label{fig:eigenface}
\end{center}
\end{figure}
%
The steps for the scheme mentioned above are given below:\\
Eigenfaces are calculated as in \cite{turk91}:\\
\underline{\textsl{Step - 1}}: Obtain face images $I_1, I_2, ... ,I_M$ (for training), each of dimension say $N\times N$.\\
\underline{\textsl{Step - 2}}: Represent every image $I_i$ as a vector $\Gamma_i$ (of dimension $N^2\times1$).\\
\underline{\textsl{Step - 3}}: Compute the average or mean face vector $\Psi$:
\begin{equation}
\Psi=\frac{1}{M}\sum_{i=1}^M \Gamma_i
\end{equation}		
\underline{\textsl{Step - 4}}: Subtract the mean face from each face vector $\Gamma_i$:
\begin{equation}
\Phi_i=\Gamma_i-\Psi
\end{equation}			
\underline{\textsl{Step - 5}}: The covariance matrix $C$ given by:
\begin{equation}
C=\frac{1}{M}\sum_{n=1}^M {\Phi_n}{\Phi_n}^T=AA^T (N^2\times N^2)
\end{equation}
	  where $A=[\phi_1,\phi_2, ... ,\phi_M]$ ($N^2\times M$ matrix) is very large, compute $A^T A$ instead ($M\times M$ matrix, $M \ll N$).\\
Now $AA^T=U\Sigma{\Sigma^TU^T}$ and $A^TA=V{\Sigma^T}\Sigma V^T$.\\
\underline{\textsl{Step - 6}}: Compute the singular vectors $v_i$ of $A^TA$. Using the equation $u_i=Av_i$ [14], singular vectors $u_i$ of $AA^T$ are obtained.\\
\underline{\textsl{Step - 7}}: Depending on energy content of the components, keep only $K$ singular vectors or eigenvectors corresponding to the $K$ largest singular values. These $K$ singular vectors are the eigenfaces corresponding to the set of $M$ face images.\\
\underline{\textsl{Step - 8}}: Normalize these $K$ eigenvectors so that $u_i=\frac{u_i}{\parallel {u_i} \parallel}_2$.\\
\underline{\textsl{Step - 9}}: Given a test image $ \Gamma $ , subtract the mean image $\Psi$ : $\Phi=\Gamma-\Psi$.\\ Project the normalized test image onto the face space
and obtain weight vector, $\Omega=U^T \Phi$, where $\Omega=[\omega_1 \omega_2 \dots \omega_k]$.\\
\underline{\textsl{Step - 10}}: Compute reconstructed image ${\Hat\Phi}=\sum_{i=1}^k {w_iu_i}$,  $(w_i={u_i}^T\Phi)$ .\\
\underline{\textsl{Step-11}}: Compute error, $e=\|\Phi-{\Hat\Phi}\|$.\\
\underline{\textsl{Step-12}}: Classify test image as belonging to face class image for which the error norm $e$ is minimum.\\
PERcentage eye CLOSure (PERCLOS) is defined as the proportion of time for which the eyelid remains closed more than 70-80\% within a predefined time period. Level of drowsiness can be judged based on the PERCLOS threshold value \cite{dinges}. Let $N_a$ be the number of eye frames belong to the open or attentive category out of $N_m$ number of eye frames captured in a minute. Hence $(N_m-N_a)$ is the number of eye frames belonging to the inattentive category. Then PERCLOS value per minute is
\begin{equation}
PERCLOS=\frac{N_m-N_a}{N_m}\times 100\%
\end{equation}
However, the PERCLOS value computed per minute is not the correct measurement of drowsiness. Literature suggests that 20 minutes bout to bout measurement is a prominent indicator of drowsiness \cite{dinges},\cite{dinges2004}. So a trade-off between execution time and accuracy of PERCLOS is requird. PERCLOS measured within a 3 minute time interval is found to indicate the level of drowsiness reasonably well. A threshold over the PERCLOS value can be used to indicate level of drowsiness \cite{dinges2004}. A comparison of eigenspace preparation using QR and QL and divide and conquer algorithm is shown in Table \ref{tab:extim}.
\begin{table}
\begin{center}
\caption{Comparison of Execution Time (sec) of QL and QR and Divide and Conquer Algorithm}
\begin{tabular}{|l|c|r|p{5cm}|}
\hline
Order & QR and QL Algorithm & Divide and Conquer Algorithm \\
\hline
50	& 0.046875	& 0 \\
100	& 0.046875	& 0.03125 \\
200	& 0.15625	& 0.03125\\
500	& 1.671875	& 0.046875\\
1000 & 9.421875	& 0.09375\\
\hline
\end{tabular}
\end{center}
\label{tab:extim}
\end{table}
\section{Conclusion and Future Scope }
In this article different symmetric SVD algorithms have been evaluated for real-time preparation of eigen space for face and eye tracking to assess the level of alertness of human driver. Different symmetric SVD algorithms have been discussed with their complexity, advantages and disadvantages. Issues related to fast Divide and Conquer SVD algorithm have been discussed. A case study of real-time face and eye tracking has been illustrated with comparison of execution time for eigen space preparation with QR and QL and Divide and Conquer SVD algorithms. The scheme has been implemented in a workstation with Intel quad processor (2.5 GHz) and 3 GB DDR RAM.\\
The present work has been undertaken with an aim to implement the algorithm in an embedded platform where processing speed and other constraints exist. Implementation of fixed-point fast SVD for face and eye tracking could be useful for fixed-point DSPs and other fixed-point embedded platforms like ARM. To convert floating-point algorithms into fixed-point,it is necessary to first estimate the range of the required floating-point variables followed by choice of Q-format. Then the floating-point variables and floating-point arithmetic equations are converted to fixed-point format \cite{zoran2007},\cite{kim98}. The fixed-point algorithm is then tested for error and the algorithm is optimized for our purpose. This fast and fixed-point SVD algorithm could be used for other embedded signal and image processing applications.Development of other fixed-point signal processing algorithms like Wavelet Transform, Hidden Markov Model could be few future extension.

%
%
%
%
%
\end{document}